\documentclass{roadef}
\usepackage{amsmath,amssymb}

\usepackage{tikz}

\usetikzlibrary{arrows, shapes}

\usetikzlibrary{arrows.meta}
\definecolor{darkred}{rgb}{0.8,0,0}
\definecolor{darkgreen}{rgb}{0,0.6,0}
\definecolor{darkblue}{rgb}{0,0,0.8}
\definecolor{darkgrey}{rgb}{0.66,0.66,0.66}
\definecolor{dauphineblue}{rgb}{0.19,0.267,0.5176}

\usepackage[font=small,labelfont=bf]{caption} 
\newcommand{\sRSP}{$1$-S-RSP}

\usepackage[colorlinks=true,bookmarks=false,citecolor=blue,urlcolor=blue,linkcolor=blue]{hyperref}

\begin{document}

\title{Resilient and Survivable Ring Star Problems}

\def\shorttitle{Resilient and Survivable Ring Star Problems}

\author{Julien Khamphousone\inst{1}, Fabian Casta\~no\inst{2}, Andr\'e Rossi\inst{1}, Sonia Toubaline\inst{1}}

\institute{
Universit\'e Paris Dauphine, PSL, France \\
\email{\{julien.khamphousone,andre.rossi,sonia.toubaline\}@dauphine.psl.eu}
\and
Frubana, Bogota, Colombia \\
\email{fabian.castano@javerianacali.edu.co}
}

\maketitle

\thispagestyle{empty}


\begin{abstract}
~In this paper, we consider both the Resilient Ring Star Problem, in which a solution should be easy to fix when a single hub fails, and the Survivable Ring Star Problem, in which a solution guarantees that a Ring Star topology is available at no cost when a single hub fails. An ILP formulation is proposed for both problems, as well as a Benders decomposition. The solution provided by both problems are also compared in order to determine which problem returns the most appropriate solutions, when the failure rate varies.
\end{abstract}

\section{Introduction}
    
\vspace*{-.1cm}\begin{minipage}{.65\textwidth}
Optimization problems involving the design of a network in a tributary or backbone architecture arise frequently in the context of transportation, telecommunication and facility location problems among many others (see Klincewicz\cite{Klincewicz98}). We consider the Ring Star network design, where a complete mixed graph with both arcs from and to every node is given, as well as edges between any pair of different nodes and a special node called the depot. The {\sc Ring Star Problem (RSP)} consists in selecting a subset of nodes including the depot, called hubs, and link them with a cycle to form the ring. Each terminal (\textit{i.e.}, non-hub) node is then connected to exactly one hub in the cycle, which is the star topology part, see Fig.~\ref{fig:small_rsp_instance} for an example. The aim of {\sc RSP} is to minimize the sum of three costs corresponding to (i) selecting the subsets of hubs, (ii) linking the ring and (iii) connecting the terminals to the ring.
\end{minipage}\hspace{.5cm}%
\begin{minipage}[h]{0.3\textwidth}

\parbox[c]{.35\textwidth}{\scalebox{.5}{%
 \begin{tikzpicture}

\begin{scope}[every node/.style={draw, thick, minimum size = 0.5cm, inner sep = 0pt, color=darkred}]
\node[rectangle] (1) at (6.5,7.5) {};
\node[circle] (2) at (3,6) {};
\node[circle] (3) at (7.5,6) {};
\node[circle] (5) at (3.5,4.5) {};
\node[circle] (6) at (7.5,2.5) {};
\node[circle] (7) at (8.5,8.5) {};

\node[circle] (9) at (4.5,3) {};

\node[circle] (hubs-draw) at (3.5,9.25) {};
\node[rectangle] (root-draw) at (3.5,10) {};
\end{scope}

\begin{scope}[every node/.style={draw, thick, minimum size = 0.5cm, inner sep = 0pt, color=teal}]
\node[circle] (non-hubs-draw) at (3.5,8.5) {};
\end{scope}

\begin{scope}[every node/.style={color=darkred}]
\node[] (hubs) at (4.4,9.25) {\large\textbf{hubs}};
\node[] (root) at (4.35,10) {\large\textbf{root}};
\end{scope}

\begin{scope}[every node/.style={color=teal}]
\node[] (non-hubs) at (4.9,8.5) {\large\textbf{terminals}};
\end{scope}

\begin{scope}[every node/.style={draw, thick, minimum size = 0.5cm, inner sep = 0pt, color=teal}]
\node[circle] (4) at (11,9) {};
\node[circle] (8) at (4.5,1.5) {};
\end{scope}
\path [draw,thick,darkred, line width=0.6mm] (1) -- (2);
\path [draw,thick,darkred, line width=0.6mm] (1) -- (7);
\path [draw,thick,darkred, line width=0.6mm] (2) -- (5);
\path [draw,thick,darkred, line width=0.6mm] (3) -- (6);
\path [draw,thick,darkred, line width=0.6mm] (3) -- (7);
\path [draw,thick,darkred, line width=0.6mm] (5) -- (9);
\path [draw,thick,darkred, line width=0.6mm] (6) -- (9);
\path [>={Stealth[teal]}, ->,draw, teal, line width=0.6mm] (4) -- (7);
\path [>={Stealth[teal]}, ->,draw, teal, line width=0.6mm] (8) -- (9);
\end{tikzpicture}}}

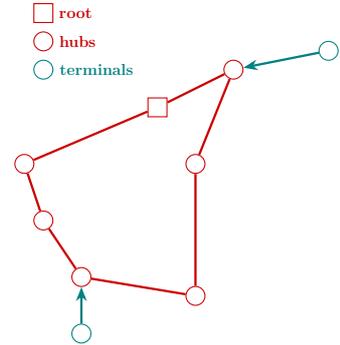
\captionof{figure}{A solution to RSP.}
\label{fig:small_rsp_instance}
\end{minipage}

{\sc RSP} has been largely studied in the literature, Labb\'e et al.\cite{labbe2004ring} proposed a Mixed Integer Programming model, strengthened with valid inequalities studied with a polyhedral analysis and solved with a Branch-and-Cut algorithm. Another exact approach that takes advantage of the fact that the depot must be in the ring can be found in Kedad-Sidhoum and Nguyen\cite{kedad2010exact}. 
    
\section{{\sc Resilient Ring Star Problem (RRSP)} and {\sc Survivable Ring Star Problem (1-S-RSP)}}

We study the {\sc Resilient Ring Star Problem (RRSP)} in the case where at most one hub can fail at any time. Each node is either \textit{certain} if it cannot fail when selected as a hub, or \textit{uncertain} if it may fail when selected as a hub (terminals are not supposed to fail). The class of nodes is a problem input{:} certain nodes are in $\widetilde{V} \subseteq V$, uncertain ones are in $V \backslash \widetilde{V}$. 

The resilient ring-star network is designed so that when an uncertain hub fails, two corrective operations occur. A ring correction consists in restoring the ring by adding a backup edge that joins the two neighbors of the failing hub, and a star repair operation aims to connect the terminals that were originally connected to the failing hub to another hub with backup arcs. Fig.~2 illustrates a solution of {\sc RRSP} and Fig.~3 presents the final ring-star structure resulting from the failure of the top-right hub. The objective of {\sc RRSP} is to minimize the usual RSP cost plus the maximum cost incurred by the corrective operations that occur when a single hub fails. The regular edges that define a solution to {\sc RSP} are expected to be used for a time period whose duration is known (typically one year). By contrast, backup edges are used to cope with the failure of a single hub, which incurs a cost that depends on the duration of the failure, hence backup costs are expressed in monetary unit per unit of time. The decision maker has to provide the input parameter $F$, which is the total amount of time during which at most one hub may be down during the time period. 

The proposed solution to {\sc RRSP} minimizes the total cost incurred by the deployment of the solution, plus the repairing operations that can occur during the time period, provided that the total failing time does not exceed $F$. A similar use of such a parameter can be found in Bertsimas and Sim\cite{bertsimas2004price}. This parameter expresses the risk protection for the solution that the decision maker wishes to achieve. Low values of $F$ favor a solution to RSP with low regular cost (RRSP reduces to RSP when $F=0$), whereas high values for $F$ lead to minimize the maximum failing cost in the RSP solution. When $F$ is large we might want to consider the survivable variant of the Ring Star Problem (\sRSP). In this variant, the backup edges and backup arcs are pre-built and always paid in advance such that whenever an uncertain hub fails, all backup arcs and edges are available for free (see Fig. 4).
    
\noindent\begin{tabular}{ccc}
 	\parbox[c]{.3\textwidth}{\scalebox{.5}{
 	\begin{tikzpicture}
 	
 	\begin{scope}[every node/.style={draw, thick, minimum size = 0.4cm, inner sep = 0pt, color=darkred}]
 	\node[rectangle] (1) at (6.5,7.5) {};
 	\node[circle] (3) at (7.5,6) {};
 	
 	\node[circle] (5) at (3.5,4.5) {};
 	\node[circle] (6) at (7.5,2.5) {};
 	\node[circle] (7) at (8.5,8.5) {};
 	
 	\node[circle] (9) at (4.5,3) {};
 	
 	\end{scope}
 	
 	\begin{scope}[every node/.style={draw, thick, minimum size = 0.4cm, inner sep = 0pt, color=darkred, fill=darkred!60}]
 	\node[circle] (3) at (7.5,6) {};
 	
 	\node[circle] (5) at (3.5,4.5) {};
 	\node[circle] (6) at (7.5,2.5) {};
 	\node[circle] (7) at (8.5,8.5) {};
 
 	
 	\end{scope}
 	
 	\begin{scope}[every node/.style={draw, thick, minimum size = 0.4cm, inner sep = 0pt, color=teal}]
 	\node[circle] (2) at (3,6) {};
 	\node[circle] (4) at (11,9) {};
 	\node[circle] (8) at (4.5,1.5) {};
 	\end{scope}
 	\path [draw,thick,darkred, line width=0.5mm] (1) -- (5);
 	\path [draw,thick,darkred, line width=0.5mm] (1) -- (7);
 	\path [draw,thick,darkred, line width=0.5mm] (3) -- (6);
 	\path [draw,thick,darkred, line width=0.5mm] (3) -- (7);
 	\path [draw,thick,darkred, line width=0.5mm] (5) -- (9);
 	\path [draw,thick,darkred, line width=0.5mm] (6) -- (9);
 	\path [draw,thick,darkred, line width=0.4mm, dashed] (1) -- (3);
 	\path [draw,thick,darkred, line width=0.4mm, dashed] (7) -- (6);
 	\path [draw,thick,darkred, line width=0.4mm, dashed] (3) -- (9);
 	\path [draw,thick,darkred, line width=0.4mm, dashed] (9) -- (1);
 	\path [>={Stealth[teal]}, ->,draw, teal, line width=0.5mm] (4) -- (7);
 	\path [>={Stealth[teal]}, ->,draw, teal, line width=0.5mm] (2) -- (5);
 	\path [>={Stealth[teal]}, ->,draw, teal, line width=0.5mm] (8) -- (9);
 	\path [>={Stealth[teal]}, ->,draw, teal, line width=0.4mm, dashed] (4) -- (1);
 	\path [>={Stealth[teal]}, ->,draw, teal, line width=0.4mm, dashed] (2) -- (9);
 	
 	\begin{scope}[every node/.style={thick, minimum size = 0.5cm, inner sep = 0pt, color=darkred}]
 	
 	\node[circle] (edges-draw) at (2.5,10.75) {};
 	\node[circle] (backupedges-draw) at (2.5,10) {};
 	\node[circle] (arcs-draw) at (2.5,9.25) {};
 	\node[circle] (backuparcs-draw) at (2.5,8.5) {};
 	\end{scope}
 	\path [draw,thick,darkred, line width=0.4mm] (edges-draw) -- (3.5,10.75);
 	\path [draw,thick,darkred, line width=0.4mm, dashed] (backupedges-draw) -- (3.5,10);
 	\path [>={Stealth[teal]}, ->,draw,thick,teal, line width=0.4mm] (arcs-draw) -- (3.5, 9.25);
 	\path [>={Stealth[teal]}, ->,draw,thick,teal, line width=0.4mm, dashed] (backuparcs-draw) -- (3.5, 8.5);
 	
 	\begin{scope}[every node/.style={color=darkred}]
 	\node[] (edges) at (4.75,10.75) {\large\textbf{ring edges}};
 	\node[] (root) at (5.05,10) {\large\textbf{backup edges}};
 	\end{scope}
 	
 	\begin{scope}[every node/.style={color=teal}]
 	\node[] (non-hubs) at (4.65,9.25) {\large\textbf{star arcs}};
 	\node[] (non-hubs) at (4.95,8.5) {\large\textbf{backup arcs}};
 	\end{scope}
 	\end{tikzpicture}}}
 	 & 
	\parbox[c]{.3\textwidth}{\scalebox{.5}{
 	 \begin{tikzpicture}
 	 	
 	\begin{scope}[every node/.style={draw, thick, minimum size = 0.4cm, inner sep = 0pt, color=darkred}]
 	\node[rectangle] (1) at (6.5,7.5) {};
 	\node[circle] (3) at (7.5,6) {};
 	
 	\node[circle] (5) at (3.5,4.5) {};
 	\node[circle] (6) at (7.5,2.5) {};
 	
 	\node[circle] (9) at (4.5,3) {};
 	
 	\end{scope}
 	
 	\begin{scope}[every node/.style={draw, thick, minimum size = 0.4cm, inner sep = 0pt, color=darkred, fill=darkred!60}]
 	\node[circle] (3) at (7.5,6) {};
 	
 	\node[circle] (5) at (3.5,4.5) {};
 	\node[circle] (6) at (7.5,2.5) {};
 	\end{scope}
 	
 	\begin{scope}[every node/.style={draw, thick, minimum size = 0.4cm, inner sep = 0pt, color=darkgrey, fill=darkgrey!20}]
 	\node[circle] (7) at (8.5,8.5) {};
 	
 	\end{scope}
 	\begin{scope}[every node/.style={draw, thick, minimum size = 0.4cm, inner sep = 0pt, color=teal}]
 	\node[circle] (2) at (3,6) {};
 	\node[circle] (4) at (11,9) {};
 	\node[circle] (8) at (4.5,1.5) {};
 	\end{scope}
 	\path [draw,thick,darkred, line width=0.5mm] (1) -- (5);
 	\path [draw,thick,red, line width=0.6mm] (1) -- (3);
 	\path [draw,thick, darkgrey, line width=0.4mm, dashed] (1) -- (7);
 	\path [draw,thick, darkgrey, line width=0.4mm, dashed] (7) -- (3);
 	\path [draw,thick,darkred, line width=0.5mm] (3) -- (6);
 	\path [draw,thick,darkred, line width=0.5mm] (5) -- (9);
 	\path [draw,thick,darkred, line width=0.5mm] (6) -- (9);
 	\path [>={Stealth[teal]}, ->,draw, teal, line width=0.6mm] (4) -- (1);
 	\path [>={Stealth[teal]}, ->,draw, teal, line width=0.5mm] (2) -- (5);
 	\path [>={Stealth[teal]}, ->,draw, teal, line width=0.5mm] (8) -- (9);
 	\begin{scope}[every node/.style={draw, thick, minimum size = 0.4cm, inner sep = 0pt, color=darkgrey}]
 	\node[circle] (failing-hub) at (3.1,10.05) {};
 	\end{scope}
 	
 	\begin{scope}[every node/.style={thick, minimum size = 0.5cm, inner sep = 0pt, color=darkred}]
 	
 	\node[circle] (edges-draw) at (2.5,10) {};
 	\node[circle] (backupedges-draw) at (2.5,10.75) {};
 	\node[circle] (arcs-draw) at (2.5,9.25) {};
 	\node[circle] (backuparcs-draw) at (2.5,8.5) {};
 	\end{scope}
 	\path [draw,darkgrey, line width=0.4mm, dashed] (backupedges-draw) -- (3.5,10.75);
	
 	\begin{scope}[every node/.style={color=darkgrey}]
 	\node[] (root) at (5.3,10.75) {\large\textbf{unavailable edges}};
 	\end{scope}
 	\begin{scope}[every node/.style={color=darkgrey}]
 	\node[] (hubs) at (4.75,10) {\large\textbf{failing hub}};
 	\end{scope}
 	\begin{scope}[every node/.style={color=darkred}]
 	\end{scope}
 	\end{tikzpicture}}}
 	&
	\parbox[c]{.3\textwidth}{\scalebox{.5}{%
 	\begin{tikzpicture}
 	
 	\begin{scope}[every node/.style={draw, thick, minimum size = 0.4cm, inner sep = 0pt, color=darkred}]
 	\node[rectangle] (1) at (6.5,7.5) {};
 	\node[circle] (3) at (7.5,6) {};
 	
 	\node[circle] (5) at (3.5,4.5) {};
 	\node[circle] (6) at (7.5,2.5) {};
 	\node[circle] (7) at (8.5,8.5) {};
 	
 	\node[circle] (9) at (4.5,3) {};
 	
 	\end{scope}
 	
 	\begin{scope}[every node/.style={draw, thick, minimum size = 0.4cm, inner sep = 0pt, color=darkred, fill=darkred!60}]
 	\node[circle] (3) at (7.5,6) {};
 	
 	\node[circle] (5) at (3.5,4.5) {};
 	\node[circle] (6) at (7.5,2.5) {};
 	\node[circle] (7) at (8.5,8.5) {};
 	
 	\end{scope}
 	
 	\begin{scope}[every node/.style={draw, thick, minimum size = 0.4cm, inner sep = 0pt, color=teal}]
 	\node[circle] (2) at (3,6) {};
 	\node[circle] (4) at (11,9) {};
 	\node[circle] (8) at (4.5,1.5) {};
 	\end{scope}
 	\path [draw,thick,darkred, line width=0.5mm] (1) -- (5);
 	\path [draw,thick,darkred, line width=0.5mm] (1) -- (7);
 	\path [draw,thick,darkred, line width=0.5mm] (3) -- (6);
 	\path [draw,thick,darkred, line width=0.5mm] (3) -- (7);
 	\path [draw,thick,darkred, line width=0.5mm] (5) -- (9);
 	\path [draw,thick,darkred, line width=0.5mm] (6) -- (9);
 	\path [draw,thick,dauphineblue, line width=0.4mm] (1) -- (3);
 	\path [draw,thick,dauphineblue, line width=0.4mm] (7) -- (6);
 	\path [draw,thick,dauphineblue, line width=0.4mm] (3) -- (9);
 	\path [draw,thick,dauphineblue, line width=0.4mm] (9) -- (1);

 	\path [>={Stealth[teal]}, ->,draw, teal, line width=0.5mm] (4) -- (7);
 	\path [>={Stealth[teal]}, ->,draw, teal, line width=0.5mm] (2) -- (5);
 	\path [>={Stealth[teal]}, ->,draw, teal, line width=0.5mm] (8) -- (9);
 	\path [>={Stealth[dauphineblue]}, ->,draw, dauphineblue, line width=0.4mm] (4) -- (1);
 	\path [>={Stealth[dauphineblue]}, ->,draw, dauphineblue, line width=0.4mm] (2) -- (9);
 
 	\begin{scope}[every node/.style={thick, minimum size = 0.5cm, inner sep = 0pt, color=darkred}]
 	
 	\node[circle] (edges-draw) at (2.5,10.75) {};
 	\node[circle] (backupedges-draw) at (2.5,10) {};
 	\node[circle] (arcs-draw) at (2.5,9.25) {};
 	\node[circle] (backuparcs-draw) at (2.5,8.5) {};
 	\end{scope}
 	\path [draw,thick,darkred, line width=0.4mm] (edges-draw) -- (3.5,10.75);
 	\path [draw,thick,dauphineblue, line width=0.4mm] (backupedges-draw) -- (3.5,10);
 	\path [>={Stealth[teal]}, ->,draw,thick,teal, line width=0.4mm] (arcs-draw) -- (3.5, 9.25);
 	\path [>={Stealth[dauphineblue]}, ->, draw,thick,dauphineblue, line width=0.4mm] (backuparcs-draw) -- (3.5, 8.5);
 	
 	\begin{scope}[every node/.style={color=dauphineblue}]
 	\node[] (root) at (5.9,10) {\large\textbf{pre-built backup edges}};
 	\end{scope}
 	
 	\begin{scope}[every node/.style={color=darkred}]
 	\node[] (edges) at (4.75,10.75) {\large\textbf{ring edges}};
 	\end{scope}
 	
 	\begin{scope}[every node/.style={color=dauphineblue}]
 	\node[] (non-hubs) at (5.8,8.5) {\large\textbf{pre-built backup arcs}};
 	\end{scope}
 	
 	\begin{scope}[every node/.style={color=teal}]
 	\node[] (non-hubs) at (4.65,9.25) {\large\textbf{star arcs}};
 	\end{scope}
 	
 	\end{tikzpicture}}} 	\\[6.5eM]
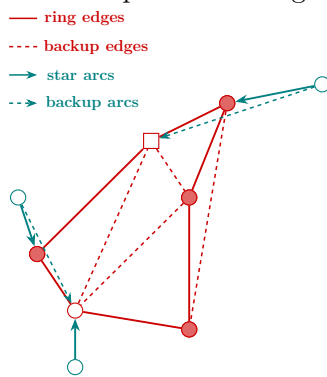
 	\parbox[t]{.3\textwidth}{\noindent\textbf{FIG. 2} -- An illustration of a solution for the {\sc RRSP} on the example of Fig. 1. where the depot and bottom left node are certain}
 	&
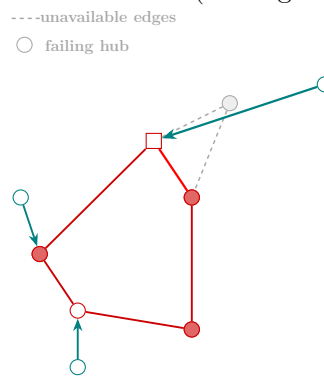
 	\parbox[t]{.3\textwidth}{\noindent\textbf{FIG. 3} -- The RSP solution that results if the gray hub fails}
 	&
 	\parbox[t]{.3\textwidth}{\noindent\textbf{FIG. 4} -- An illustration of a solution for the {\sRSP} on the example of Fig. 1. where the depot and bottom left node are certain}\\
 \end{tabular}
 
We have formulated and implemented a Mixed Integer Linear Programming model for both RRSP and \sRSP~with the set $\widetilde V$ that contains the uncertain nodes. We compare the results of this formulation and a Benders decomposition, and analyze the impact of parameter $F$ in the comparative cost of the optimal solutions to both problems.

\bibliographystyle{plain}

\end{document}